\documentclass{article}

\usepackage{amsmath}
\usepackage{amsfonts}
\usepackage{amssymb}
\usepackage{enumerate} 

\newtheorem{theorem}{Theorem}

\newtheorem{lemma}[theorem]{Lemma}

\begin{document}

\title{Almost Repdigits in $ k-$generalized Lucas Sequences}
\author{Alaa ALTASSAN and Murat ALAN\\
King Abdulaziz University\\
 Department of Mathematics, Jeddah, 21589, Saudi Arabia \\
Yildiz Technical University\\
Department of Mathematics, Istanbul, 34210, Turkey\\
}

\maketitle

\abstract{Let $ k \geq 2 $ and  $ ( L_{n}^{(k)} )_{n \geq 2-k} $ be the $k-$generalized Lucas sequence with initial condition $ L_{2-k}^{(k)} = \cdots = L_{-1}^{(k)}=0 ,$ $ L_{0}^{(k,}=2,$ $ L_{1}^{(k)}=1$ and each term afterwards is the sum of the $ k $ preceding terms. A positive integer is an almost repdigit if its digits are all equal except for at most one digit. In this paper, we work on the problem of determining all terms of $k-$generalized Lucas sequences which are almost repdigits. In particular, we find all $k-$generalized  Lucas numbers which are powers of $10$ as a special case of almost repdigits.}






\maketitle

\section{Introduction}

Let $ k \geq 2 $ be an integer. The $ k- $generalized Lucas sequence or, for simplicity, the $k-$Lucas sequence is a sequence given by the recurrence relation
$$ L_{n}^{(k)}=L_{n-1}^{(k)} + \cdots + L_{n-k}^{(k)} \quad \text{for all} \quad n \geq 2 ,$$
with the initial values $ L_{i}^{(k)}=0 $ for $ i=2-k, \ldots, -1 ,$ $ L_{0}^{(k)}=2$  and $ L_{1}^{(k)}=1.$ We call the terms of this sequence $ k-$Lucas numbers for simplicity.  For $ k=2, $ this sequence is the classical Lucas sequence, and hence the $ k-$Lucas sequence is a generalization of the Lucas sequence from binary recurrence sequence to the order $ k $ recurrence sequence.

Recall that, a positive integer whose  all digits are equal is called a repdigit. In recent years, many researches have been performed to find  all terms of some sequences related to repdigits, see for example \cite{Alahmadi,Bed,BBG21,BL13,BL14,BL,Coufal,Luca2000,Marques15,24,25R,ZK22}. In this study, we search the numbers similar to the repdigits in $ k-$Lucas numbers.

In \cite{Gica} and \cite{Kih}, all square and perfect power positive integers whose digits are all equal except for one digit have been examined without giving a specific name. We call a positive integer an almost repdigit if its all digits are equal except for at most one digit. This numbers can be written of the form
$$  a \left( \dfrac{10^{d_{1}} -1 }{9} \right)+(b-a)10^{d_{2}}, \quad 0 \leq  d_{2} <  d_{1} \quad \text{and} \quad 0 \leq a, b \leq 9.$$
Usual repdigits and the numbers of the form $ b10^{d_2} $ are two particular cases of almost repdigits which corresponds to the cases $ a=0 $ and $ b=a ,$ respectively. Thus, almost repdigits are a generalization of repdigits.

Recently, in \cite{Altassan}, the authors found all $ k-$generalized Fibonacci numbers that are almost repdigits. In this paper, we continue to search almost repdigits by taking into account $ k-$generalized Lucas numbers as an analogue of the study in \cite{Altassan}. In other words, we consider the Diophantine equation
\begin{equation}
L_n^{(k)} =a \left( \dfrac{10^{d_{1}}-1}{9} \right)+(b-a)10^{d_{2}} , \quad  0 \leq  d_{2} <  d_{1} \quad \text{and} \quad 0 \leq a, b \leq 9
\label{Lnab}
\end{equation}
in non negative integers $d_{1},d_{2},a$ and $ b $ and we prove the following theorem.
\begin{theorem}
\label{main2}
The Diophantine equation \eqref{Lnab} has solutions only in the cases $ L_{11}^{(2)}=199 ,$ $ L_{12}^{(2)}=322 ,$ $ L_{8}^{(3)}=118 ,$ $ L_{10}^{(3)}=399 ,$ $ L_{10}^{(7)}=755 $ and $ L_{10}^{(9)}=766 ,$ when $ L_n^{(k)}$ has at least three digits.
\end{theorem}

Since, the numbers having at most two digits are trivially almost repdigits, we state the above theorem only for $ L_n^{(k)}$  which are consist of at least three digits. Thus, from now on we take $ d_1 \geq 3$ and hence $ n>5.$

The proof of the above theorem, mainly depends on two effective methods, that is, the linear forms in logarithms of algebraic numbers due to Matveev \cite{Matveev} as well as reduction algorithm due to Dujella and Peth\H o  \cite{DP}, which is in fact originally introduced by Baker and Davenport in \cite{Baker-Davenport}. We give some details of these methods in the next section whereas in the third section, we give the main properties of $k-$Lucas sequences that we will need later. We devoted the forth section to the proof of Theorem \ref{Lnab}. It is also worth to note that, we implemented the software Maple for all calculations and computations in the proof the Theorem.

\section{The Tools}

Let $\theta$ be an algebraic number, and let
\[
c_0x^d+c_1x^{d-1}+\cdots+c_d=c_0\prod_{i=1}^{d}(x-\theta^{(i)})
\]
be its minimal polynomial over $ \mathbb{Z} ,$ with degree $d,$ where the $c_i$'s are relatively prime integers with $c_0>0,$ and the $\theta^{(i)}$'s are conjugates of $\theta$.

The logarithmic height of $\theta$ is defined by
$$
h(\theta)=\dfrac{1}{d}\left(\log c_0+\sum_{i=1}^{d}\log\left(\max\{ \lvert \theta^{(i)} \rvert ,1\}\right)\right).
$$
If $\theta= r/s$ is a rational number with relatively prime integers $ r $ and $ s $ and $s>0,$ then  $h(r/s)=\log \max \{|r|,s\}$. The following properties are very useful in calculation of a logarithmic height :
\begin{itemize}
\item[$ \bullet $] $h(\theta_1 \pm \theta_2)\leq h(\theta_1) + h(\theta_2)+\log 2$.
\item[$ \bullet $] $h(\theta_1 \theta_2^{\pm 1})\leq h(\theta_1)+h(\theta_2)$.
\item[$ \bullet $] $h(\theta^{s})=|s|h(\theta),$ $ s \in \mathbb{Z} $.
\end{itemize}
\begin{theorem}[Matveev's Theorem]
\label{Matveev}
Assume that $\alpha_1, \ldots, \alpha_t$ are positive real algebraic numbers in a real algebraic number field $\mathbb{K}$ of degree $ d_\mathbb{K} $ and let $b_1, \ldots, b_t$ be rational integers, such that 
\[
\Lambda:=\alpha_1^{b_1}\cdots\alpha_t^{b_t}-1,
\]
is not zero. Then
$$
\lvert \Lambda\rvert >\exp\left( K(t) d_\mathbb{K}^2(1+\log d_\mathbb{K})(1+\log B)A_1\cdots A_t\right),
$$
where
\[
K(t):= -1.4\cdot 30^{t+3}\cdot t^{4.5}  \quad \text{and} \quad    B\geq \max\{|b_1|,\ldots,|b_t|\},
\]
and
\[
A_i\geq \max\{d_\mathbb{K}h(\alpha_i),\lvert \log \alpha_i \rvert, 0.16\}, \quad \text{for all} \quad  i=1,\ldots,t.
\] 
\end{theorem}
For a real number $ \theta, $ we put $ \lvert \lvert \theta \rvert \rvert =\min\{ \lvert \theta -n \rvert :  n \in \mathbb{Z} \}, $ which represents the distance from $ \theta $ to the nearest integer. Now, we cite the following lemma which we will use it to reduce some upper bounds on the variables.

\begin{lemma} \label{reduction}\cite[Lemma 1]{13}      
Let $M$ be a positive integer, and let $p/q$ be a convergent of the continued fraction of the irrational $\gamma$ such that $q>6M.$ Let $A,B,\mu$ be some real numbers with $A>0$ and $B>1$. If $\epsilon:= \lvert \lvert \mu q \rvert \rvert -M \lvert \lvert \gamma q \rvert \rvert >0$, then there is no solution to the inequality
\[
0< \lvert u\gamma-v+\mu \rvert <AB^{-w},
\]
in positive integers $u,v$ and $w$ with
\[
u\leq M \quad\text{and}\quad w\geq \frac{\log(Aq/\epsilon)}{\log B}.
\]
\end{lemma}

We cite the following lemma from \cite[Lemma 7]{Guzman}.
\begin{lemma}
\label{Guzman}
Let $ m \geq 1 $ and $ T>(4m^2)^m .$ Then we have 
$$
\dfrac{x}{(\log x)^m} < T \Rightarrow x < 2^mT(\log(T))^m.
$$
\end{lemma}

\section{Properties of  $k-$Fibonacci and $k-$Lucas Numbers }

The characteristic polynomial of $ k- $Lucas numbers is
$$\Psi_k(x)=x^k-x^{k-1}- \cdots -x-1$$ 
which is an irreducible polynomial over $ \mathbb{Q}[x] .$ The polynomial $ \Psi_k(x) $ has exactly one real distinguished root $ \alpha(k) $ outside the unit circle \cite{22,23,26}.  The other roots of $ \Psi_k(x) $ are strictly inside the unit circle \cite{23}. This root $ \alpha(k) ,$ say  $ \alpha $ for simplicity, placed in the interval
$$ 2(1-2^{-k})< \alpha < 2  \quad \text{for all} \quad k \geq 2. $$

Let
\begin{equation}\label{fk}
f_k (x) =  \dfrac{x - 1}{2 + (k+1)(x - 2)}.
\end{equation}
It is known that the inequalities
\begin{equation}\label{fkprop}
1/2 <f_k (\alpha) < 3/4  \quad \text{and} \quad \left \lvert f_k( \alpha_{i} ) \right\rvert <1, \quad 2 \leq i \leq k
\end{equation}
are hold,  where $ \alpha:=\alpha_{1}, \cdots , \alpha_{k}  $ are all the roots of $ \Psi_k(x) $  
\cite[Lemma 2]{13}. In particular, we deduce that $  f_k(\alpha) $ is not an algebraic integer. In the same Lemma, it is also proved that
\begin{equation}\label{hfk}
h( f_k (\alpha) ) < 3 \log(k) \quad \text{holds} \quad \forall k \geq 2.
\end{equation}

In \cite{BL14}, Bravo and Luca showed that
\begin{equation}\label{kL}
L_n^{(k)} = \sum_{i=1}^k  (2 \alpha_{i}-1)   f_k (\alpha) (\alpha_{i})^{n-1} \quad \text{and} \quad \left \lvert L_n^{(k)} - (2 \alpha-1)  f_k (\alpha) \right \rvert < 3/2
\end{equation}
for all $ k\geq 2. $
As in the classical $ k=2 $ case, we have the similar bounds as
\begin{equation}
\label{L1}
\alpha^{n-1} \leq L_n^{(k)} \leq 2\alpha^{n} 
\end{equation}
for all $n\geq 1$ and $ k \geq 2 $  \cite{BL14}.

\section{Proof of Theorem \ref{main2}}

Assume that Equation \eqref{Lnab} holds. By combining the inequality
$$ 10^{{d_1}-2} <  a \left( \dfrac{10^{d_{1}} -1}{9} \right)+(b-a)10^{d_{2}}  \leq 2 \cdot 10^{d_1} ,$$
and \eqref{L1}, we get that
\begin{equation}
\label{d1nL}
{d_1}  <  \dfrac{\log 2 }{\log 10} (n+1) +2 < 0.31n + 2.31 < n-1
\end{equation}
and
\begin{equation}
\label{d1n2L}
0.2n-0.6 < \dfrac{\log ((1+\sqrt(5))/2)}{\log 10} (n-1)-\dfrac{\log 2}{\log 10} <   {d_1}  
\end{equation} 
for all $ n>5. $

First, we assume that $ a \neq 0 .$ We examine the case $ a=0 $ in the end of this section.

\subsection{The Case $ n < k+1 $ and Almost Repdigits of the Form $3 \cdot 2^n$}

Assume that $ n \leq k. $ In this case, $ L_n^{(k)} = 3 \cdot 2^{n-2}$, and hence Equation \eqref{Lnab} can be written as
$$ 27 \cdot 2^{n-2}  =a \left( 10^{d_1} -1 \right)+9(b-a)10^{d_{2}} .$$
Taking modulo $ 2^{d_{2}} $ and modulo $ 2^{d_{1}} ,$ we find that $ {d_{2}} \leq 3 $ and $ {d_{1}} \leq 13. $ 

Hence, from \eqref{d1n2L}, we see that $ n<70. $ A quick calculation shows that when $ n<70 ,$  there is no almost repdigits of the form $ 3 \cdot 2^{n-2} $ with at least three digits.

So, from now on, we take $ n \geq k+1. $

\subsection{A Bound for $ n $ Depending on $ k $}

By rewriting \eqref{Lnab} as

$$ L_n^{(k)} +a/9 -(b-a)10^{d_2}=a10^{d_1}/9   ,$$
and by using \eqref{kL}, we get

$$ \left\lvert  f_k(\alpha) ( 2 \alpha -1) \alpha^{n-1}-a10^{d_1}/9 \right \rvert \leq (3/2)+\left\lvert a/9 -(b-a)10^{d_2}   \right \rvert .
$$
Therefore, we obtain
\begin{equation}
\label{lam1}
\lvert \Lambda_{1} \rvert \leq \dfrac{27/2}{10^{d_1}} + \dfrac{1}{10^{d_1}} +\dfrac{\lvert b-a \rvert (9/a)}{10^{d_1-d_2}} \leq \dfrac{87}{10^{d_1-d_2}},
\end{equation}
where 
$$  \Lambda_{1} := \alpha^{n-1} 10^{-d_1}  f_k(\alpha)( 2 \alpha -1) 9/a -1.   $$

We take 
$$ \eta_1:=\alpha, \eta_2:=10, \eta_3:= f_k(\alpha) ( 2 \alpha -1) 9/a ,$$
$$ b_1:=n-1, b_2:=-d_1 ,b_3:=1. $$ 
Note that from \eqref{hfk}, we find 
$h(\eta_3) \leq  h(9/a) + h( f_k(\alpha) ) +h(2 \alpha -1)  < 8 \log(k),$
since $ h(2 \alpha -1)< \log 3 $ \cite[page 147]{BL14}.

We have also $ \Lambda_{1} \neq 0.$ Indeed, if $ \Lambda_{1}=0 ,$ then we would get
$$ a10^{d_1}/9 = f_k(\alpha) (2 \alpha -1) \alpha^{n-1} .$$
Conjugating both sides of this relation by  any one of the automorphisms  $ \sigma_i : \alpha \rightarrow \alpha_i $ for any $ i \geq 2 $ and by taking the absolute values, we find that
$$ 100 < \lvert f_k(\alpha_i) \rvert   \lvert 2 \alpha -1 \rvert   \lvert \alpha_i \rvert^{n-1} < 3,$$
a contradiction. Thus,$ \Lambda_{1} \neq 0.$ Other calculations are doneby using similar techniques as in the $k-$Fibonacci case. So, by combining the result of Theorem \ref{Matveev} and the fact that
 $ \log{(  \Lambda_{1}) } <  \log 87 - ({d_{1}-d_{2}}) \log 10, $
we obtain
\begin{equation}
\label{d12firstL}
{d_{1}-d_{2}} < 4.8 \cdot 10^{12} \cdot k^4 (\log ^2 k ) \log {(n-1)}.
\end{equation}
By rearranging Equation \eqref{Lnab} as follows
$$ L_n^{(k)} + a/9 = a10^{d_1}/9 + (b-a)10^{d_2},$$
and using \eqref{kL}, we get

\begin{equation}
\label{lam2L}
\lvert \Lambda_{2} \rvert  \leq  \dfrac{5}{2}\dfrac{1}{f_k(\alpha) \alpha^{n-1} } \leq \dfrac{5}{2 \alpha^{n-1}}.
\end{equation}
where
$\Lambda_{2} := \alpha^{-(n-1) } 10^{d_1} f_k(\alpha)^{-1} (2 \alpha -1)^{-1}   ( (a/9)+  (b-a)10^{d_2-d_1}) -1.$
By the similar argument as above we see that $ \Lambda_{2} \neq 0. $ Let 
$$ \eta_1:= \alpha, \eta_2:= 10, \eta_3:=f_k(\alpha)^{-1}  (2 \alpha -1)^{-1} ( (a/9)+  (b-a)10^{d_2-d_1}) $$
with $ b_1:= -(n-1)$, $b_2:=d_1,$  $b_3:=1. $ All $ \eta_1, \eta_2$ and $ \eta_3 $  belong to the real number field $ \mathbb{K}=\mathbb{Q}(\alpha) $ and therefore we take $ d_\mathbb{K}=2 $, to be the degree of the number field $ \mathbb{K} .$ Using the properties of logarithmic height, and the fact that $ h(2 \alpha -1 )< \log 3 ,$ we find
\begin{align*}
h(\eta_3) & \leq h(f_k(\alpha)^{-1})+ h( (2 \alpha -1)^{-1} )    + h( (a/9)+  (b-a)10^{d_2-d_1} ) \\
& \leq 3 \log(k) + h( (2 \alpha -1) )+ h(a/9)+h(b-a)+h( 10^{d_{2}-d_{1}} ) +\log(2)\\
&\leq 3 \log(k) + \log 3 + \log(144) +\lvert  d_{2}-d_{1} \rvert \log(10)\\
& < 12 \log(k) +\lvert  d_{2}-d_{1} \rvert \log(10).
\end{align*}

By Theorem \ref{Matveev}, we get a bound for $ \log(\Lambda_2).$ Then by combining this bound with the one comes from \eqref{lam2L}, we get

$$ n-1 < 2.4 \cdot 10^{25} k^8 (\log(k))^3 (\log(n-1))^2 .    $$
Therefore, from Lemma \ref{Guzman}, we have
\begin{equation}\label{nfirstL}
n< 1.3 \cdot 10^{30} k^8 (\log(k))^5.
\end{equation}

\subsection{The Case $ k \leq 470 $}

Let
$$ \Gamma_{1} := (n-1) \log \alpha -d_1 \log 10 +\log(f_k(\alpha)(2 \alpha -1)\cdot 9/a) .$$
Then 
$$ \Lambda_{1} := \left \lvert \exp(\Gamma_{1}) -1 \right \rvert < 87/10^{d_1-d_2}.$$
We claim that ${d_1-d_2}<65.$ Suppose that $ {d_1-d_2}>3.$ Then  $ 87/10^{d_1-d_2}<1/2 $ and therefore $\lvert \Gamma_1 \rvert <174/10^{d_1-d_2} .$  So we have
\begin{equation}
\label{g1L2}
0 < \left \lvert (n-1) \dfrac{\log \alpha }{\log 10 }  -d_1 + \dfrac{\log(f_k(\alpha) \cdot (2 \alpha -1) \cdot 9/a)}{\log 10 }    \right \rvert < 174/10^{d_1-d_2} \log 10.
\end{equation}

For each $ 2 \leq k \leq 470 $, we take $ M_k := 1.3 \cdot 10^{30} k^8 (\log(k))^5 >n $ and $\tau_k = \dfrac{ \log{\alpha}}{\log 10} .$ Then, for each $ k $, we find a convergent $ p_i/q_i $ of the continued fraction of irrational $\tau_k$ such that $q_i >6M_k$

After that, we calculate
$\epsilon_{(k,a)} := ||\mu_{(k,a)}  q_{i} || -M_k || \tau_k q_{i} ||  $ for each $ a \in \{1, 2, \ldots, 9 \}, $ where
$$ 
\mu_{(k,a)}  := \dfrac{\log(f_k(\alpha)\cdot (2 \alpha -1) \cdot 9/a)}{\log 10 }.
$$
If $\epsilon_{(k,a)}< 0 ,$ then we repeat the same calculation for $ q_{i+1}. $ Except for $ (k,a)=(2,9).$   In fact $ 0.00008 <\epsilon_{(k,a)}  .$ Thus, from  Lemma \ref{reduction}, we find an upper bound on $ d_1-d_2 $  for each $ 2 \leq k \leq 470 $ such that none of them are greater than $ 61. $ So, we conclude that $ d_1-d_2 < 65 $ as we claimed.

If $ (k,a)=(2,9) ,$ then $ \tau =\mu_{(k,a)} $ and hence $\epsilon_{(k,a)}= 0 .$ So, in this case, \eqref{g1L2} is

\begin{equation*}
\label{g1L}
0 <  \left \lvert n \dfrac{\log \alpha }{\log 10 }  -d_1  \right \rvert <  174/10^{d_1-d_2} \log 10.
\end{equation*}
In other words,
\begin{equation}
0 < \left \lvert  \dfrac{\log \alpha }{\log 10 }  - \dfrac{d_1}{n} \right \rvert <  174/10^{d_1-d_2} \log 10.
\end{equation}
From, the theory of continued fractions, we see that this implies that the bound $ d_1-d_2 < 65 $ is valid.

Next, let
$$ \Gamma_{2} := -(n-1) \log \alpha +d_1 \log 10 +\log( f_k(\alpha)^{-1}\cdot (2 \alpha -1)^{-1} \cdot( (a/9)+  (b-a)10^{d_2-d_1})   ) .$$
Thus,
$$ \Lambda_{2} := \left \lvert \exp(\Gamma_2) -1 \right \rvert < 5/2 \alpha^{n-1} <1/2.$$
Hence, we get that
\begin{equation}
\label{g2L}
0 < \left \lvert (n-1) \dfrac{\log \alpha }{\log 10 }  -d_1 +\mu_{(k,{d_1-d_2},a,b)}    \right \rvert < \dfrac{5}{ \alpha^{n-1} \log 10  }.
\end{equation}
where
$$ 
\mu_{(k,{d_1-d_2},a,b)}  := - \dfrac{  \log( f_k(\alpha)^{-1}\cdot (2 \alpha -1)^{-1} \cdot ( (a/9)+  (b-a)10^{d_2-d_1})   )   }{\log 10 } .
$$

This time we calculate
$\epsilon_{(k,{d_1-d_2},a,b)} := ||\mu_{(k,{d_1-d_2},a,b)} q_{i} || -M_k || \tau_k q_{i} ||  $ for each $ d_1-d_2 \in \{1, 2, \ldots, 65 \} ,$  $ a \in \{1, 2, \ldots, 9 \} $ and $ b \in \{0, 1, \ldots, 9 \} .$

We apply Lemma \ref{reduction} to Equation \eqref{g2L}, and therefore we find an upper bound on $ n-1 $  for each $ 2 \leq k \leq 470,$ say $ n_L(k). $  For example $ n_L(3)<150, $ $ n_L(10)<147, $ $ n_L(100)<178, $ $ n_L(200)<197, $  $ n_L(300)<296,$  $ n_L(400)<396 $ and  $ n(470)<465 $ are some of these bounds. 

By writing a  short computer programme in Maple, and using the obtained bounds,  we find that $ L_{11}^{(2)}=199 ,$ $ L_{12}^{(2)}=322 ,$ $ L_{8}^{(3)}=118 ,$ $ L_{10}^{(3)}=399 ,$ $ L_{10}^{(7)}=755 $ and $ L_{10}^{(9)}=766 ,$ are the only $ k- $Lucas numbers which are almost repdigits with at least three digits, as we claimed in Theorem \ref{main2}.
Now, we turn our focus to the case $ k>470. $

\subsection{The Case $ k > 470 $}

We use the following lemma \cite[Lemma 2.6]{24}.
\begin{lemma}
If $ n<2^{k/2} ,$ then the following estimates hold:
$$  L_n^{(k)}=3 \cdot 2^{n-2}(1+\zeta (n,k)), \quad \text{where} \quad \lvert \zeta (n,k) \rvert  < \dfrac{1}{2^{k/2}} .$$
\end{lemma}

For $ k>470 ,$ the inequality 
$$ n< 1.3 \cdot 10^{30} k^8 (\log(k))^5 < 2^{k/2}, $$ 
holds and hence from the above Lemma, we have that 
\begin{equation} \label{41L}
\left \lvert 3 \cdot 2^{n-2} -  L_n^{(k)} \right \rvert <3 \cdot \dfrac{2^{n-2}}{2^{k/2}}  .
\end{equation}

Now, we turn back to  \eqref{Lnab} one more time to rewrite it as
\begin{equation} \label{42L}
\left \lvert  L_n^{(k)} - (a/9)10^{d_1}  \right \rvert < (a/9)+\lvert b-a \rvert 10^{d_2}. 
\end{equation}

Thus, by combining \eqref{41L} and \eqref{42L}, we get 
$$\left \lvert 3 \cdot 2^{n-2} - (a/9)10^{d_1}  \right \rvert < 3 \cdot \dfrac{2^{n-2}}{2^{k/2}} + (a/9)+\lvert b-a \rvert 10^{d_2} .$$

Therefore, we have
\begin{equation}
\label{lam3L}
\Lambda_{3} := \left \lvert 2^{n-2} 10^{-d_1}27/a  - 1  \right \rvert < \dfrac{1}{2^\lambda} ,
\end{equation}
where $ \lambda :=\min\{ (k/2)-5, {(d_1-d_2)} \dfrac{\log(10)}{\log(2)}-8\}  .$

Let  $ \eta_1:=2,$ $\eta_2:=10,$ $ \eta_3:= 27/a $ and $ b_1:=n-2,$  $b_2:=-d_1 ,$ $b_3:=1. $ 
Applying Theorem \ref{Matveev} to $ \Lambda_{3} $, we get 
$ \lambda  < 2.2 \cdot 10^{12} \log n $

where we used
\begin{align*}
\log n  & <  \log ( 1.3 \cdot 10^{30} k^8 (\log(k))^5  ) \\
& < \log(4.3)+ 3 \log(10)+8 \log k +5 \log \log(k) \\
& < 50 \log k.
\end{align*}

Thus, if $ \lambda:=(k/2)-5 ,$ then we get a bound for $ k $  as
$$ k < 10^{16}. $$
If $ \lambda:={(d_1-d_2)} \dfrac{\log(10)}{\log(2)}-8  ,$ then we get 
\begin{equation} \label{d37L}
{d_1-d_2}< 3.4 \cdot 10^{13} \log k. 
\end{equation}

This bound of ${d_1-d_2} $ also leads to an upper bound of $ k. $ To do this, we rewrite \eqref{Lnab} as
\begin{equation} \label{43L}
\left \lvert  L_n^{(k)} - (a/9)10^{d_1} - (b-a)10^{d_2}  \right \rvert \leq (a/9) \leq 1.
\end{equation}
By combining \eqref{43L} with \eqref{41L}, we find that

$$
\label{lam4}
\Lambda_{4} := \left \lvert  2^{-(n-2)}10^{d_1}  ( (a/9)+  (b-a)10^{d_2-d_1})(1/3) -1 \right \rvert < \dfrac{1}{3 \cdot 2^{n-2} }      \dfrac{2}{2^{k/2}}  \leq \dfrac{2}{2^{k/2}}.
$$

We take  
$$ \eta_1:=2, \eta_2:=10, \eta_3:= (1/3)((a/9)+ (b-a)10^{d_2-d_1}) $$ and  $b_1:=-(n-2), b_2:=-d_1 ,b_3:=1. $ Then

\begin{align*}
h(\eta_3) &= h(a/9)+h(b-a)+\lvert {d_2-d_1} \rvert h(10) + h(3)+ \log 2 \\
& < \log 432 +({d_1-d_2} )\log 10.
\end{align*}

Other calculations are similar to those for $ \Lambda_{4} $ as $ \mathbb{K}=\mathbb{Q} ,$ $ d_ \mathbb{K}=1 ,$ $ B: = n>n-2.$  $ h(\eta_1) = \log 2, $ $ h(\eta_2) = \log 10 .$

Moreover $ \Lambda_{4} \neq 0 .$ Indeed, $2^{n-2}=10^{d_1} (a/9)+  (b-a)10^{d_2} $  implies that  $ a=9 $ and $ d_2=0. $ For $ d_1=3  $, clearly the equation $2^{n-2}=10^{d_1} + b-9 $ has no solution in integers. So $ d_1 >3. $ Thus, congruence consideration modulo $ 2^4 $ shows that this equation has no integer solutions for $ 0 \leq b \leq 9. $ Hence, $ \Lambda_{4} \neq 0 .$

Therefore,  Theorem \ref{Matveev} together with \eqref{lam4} give that
$$
\log2 - (k/2) \log2 > -1.4 \cdot 30^6 \cdot 3^{4.5}  (1+\log n) \log 2 \cdot \log 10  \cdot (\log 432 +({d_1-d_2} )\log 10 ) .
$$
At this point, we use the upper bound of $ d_1-d_2 $ which was given in \eqref{d37L}, and by using the estimates  $ \log 432 < \log k $ and $ \log n < 50 \log k, $
we obtain the desired upper bound for $ k $ as
\begin{equation}
k< 3\cdot 10^{31}.
\end{equation}
Thus, by \eqref{nfirstL}, we have also a bound for $ n $ as
\begin{equation}\label{n292L}
n < 1.8 \cdot 10^{291}.
\end{equation}

\subsection{Reducing the Bound on k}
We will reduce these highly large upper bounds.
Let 
\begin{equation}\label{g3L}
\Gamma_{3} := (n-2) \log 2 -d_1 \log 10 +\log(27/a).
\end{equation}

Then $ \Lambda_{3} :=\left \lvert \exp(\Gamma_{3})-1 \right \rvert < \dfrac{1}{2^\lambda} .$
We will find a feasible bound for $\lambda .$ Suppose that $ \lambda > 1 .$ Then, $ \dfrac{1}{2^\lambda} < \dfrac{1}{2}$ and hence we get that $ \lvert \Gamma_{3} \rvert < \dfrac{2}{2^\lambda} .$

In this case, we don't need to consider the case $ a=9 $ separately.

From \eqref{g3L}, we write
\begin{equation}
0< \left \lvert (n-2) \dfrac{ \log 2 }{\log 10}  - d_1 + \dfrac{ \log (27/a) }{ \log 10} \right \rvert < \dfrac{2}{ 2^\lambda \log 10 }.
\end{equation}

Let $ M:= 1.8 \cdot 10^{291} >n $ and  $\tau = \dfrac{ \log{2}}{\log 10}. $ Then, the denominator of the $ 588th $ convergent of $ \tau,$ say $q_{588},$ exceeds $ 6M.$ 

Then
$$ \epsilon_a := ||\mu_a q_{588} || -M || \tau q_{588} || > 0.029559 ,$$ for each $ a \in \{1, 2, \ldots, 9 \}, $ 
where
$$ \mu_a := { \log(27/a) }/{ \log10} .$$

Thus, by applying  Lemma \ref{reduction}, we get $ \lambda < 975 .$

Hence, if $ \lambda=k/2 -5,$ then $ k <1960.$ Assume that $ \lambda ={(d_1-d_2)} \dfrac{\log(10)}{\log(2)}-8 .$ Then
$$ {d_1-d_2} < 296 <300.$$
Let
\begin{equation}
\Gamma_{4} = \left \lvert (n-2) \log 2 - d_1 \log 10  -  \log( (a/9)+(b-a)10^{d_2-d_1 } )(1/3) \right \rvert 
\end{equation} 
Then
\begin{equation}
\Lambda_{4}:= \lvert \exp( \Gamma_{4} )-1 \rvert < \dfrac{2}{2^{k/2}}<\dfrac{1}{2}.
\end{equation}
So
\begin{equation}\label{g4L}
0<\lvert \dfrac{\Gamma_{4} }{\log 10}  \rvert < \left \lvert (n-2) \dfrac{\log 2 }{\log 10 } -d_1 + \mu_{(a,b,d_1-d_2)} \right \rvert  < \dfrac{4}{2^{k/2} {\log 10 }}
\end{equation} 
where
$$ \mu_{(a,b,d_1-d_2)} := - \dfrac{ \log{ ( (a/9)+(b-a)10^{d_2-d_1 } )(1/3) } }{{\log 10}}. $$

We take $ M:= 1.8 \cdot 10^{291} >n $ and  $\tau = \dfrac{ \log{2}}{\log 10}. $ 

This time, we take $ q_{595},$ which is the denominator of the $ 595th $ convergent of $ \tau,$ as $q_i > 6M.$

Let
$\epsilon_{(a,b,d_1-d_2)}:= ||\mu_{(a,b,d_1-d_2)} q_{595} || -M || \tau q_{595} ||  $ for each $ a \in \{1, 2, \ldots, 9 \} $, $ b \in \{1, 2, \ldots, 9 \} $ and $ d_1-d_2 \in \{1, 2, \ldots, 300 \} .$ 
We find that $0.000036< \epsilon_{(6,0,272)}  \leq  \epsilon_{(a,b,d_1-d_2)}$ for  all $a, b, {d_1-d_2 } $ except for $(a, b, {d_1-d_2 })=(9,2,1), (9,5,1), (9,5,2) $
since these three triples of $(a, b, {d_1-d_2 })$, we have that $ \mu_{(a,b,d_1-d_2)}<0. $

Then, except for the above three triples, from  Lemma \ref{reduction}, we conclude that $ k < 2000 .$ Hence, from \eqref{nfirstL},  $ n < 8.5 \cdot 10^{60} .$

If $(a, b, {d_1-d_2 })=(9,2,1),$ then

$$ \mu_{(9,2,1)} :=  -\dfrac{ \log{ ( 1-\dfrac{7}{10} ) (1/3) } }{{\log 10}} =1 \in \mathbb{Z} .$$
So, in this case we may write 

\begin{align*}
\Gamma_{4}  & = \left \lvert  (n-2) \log 2 - d_1 \log 10   -  \log 10^{-1} \right \rvert \\
& = \left \lvert (n-2) \log 2 - ( d_1 - 1) \log 10 \right \rvert,\\
\end{align*}
and hence
\begin{equation}
\label{gcf}
0<\left \lvert       \dfrac{\log 2 }{\log 10 }  - \dfrac{d_1 - 1}{n-2} \right \rvert  < \dfrac{4}{2^{k/2} {\log 10 }}.
\end{equation}

The inequality
$$ \dfrac{4}{ 2^{k/2} {\log 10 }   }  > \dfrac{1}{ 2 \cdot (n-2)^2 } $$
implies that $ k<2000. $

Assume that
$$ \dfrac{4}{ 2^{k/2} {\log 10 }   }  \leq \dfrac{1}{ 2 \cdot (n-2)^2 } .$$
Then $ \dfrac{d_1 - 1}{n-2} $ is a convergent of $ \dfrac{\log 2 }{\log 10 } ,$ say $ p_i /q_i. $ Then $ q_i <n-2 < 1.8 \cdot 10^{291}$ implies $ i<588 $  and   $ \max{a_i}=5393 .$ So, from the properties of continued fractions, see \cite[Theorem 1.1.(iv)]{Hen},
$$ 2^{k/2}  < \dfrac{4 \cdot 5395 \cdot 1.8 \cdot 10^{291} }{\log 10} < 1.7 \cdot 10^{295} < 2^{981}.$$
Thus, the upper bound $ k<2000 $ is valid in this case also.

If $(a, b, {d_1-d_2 } )=(9,5,1)$ then

$$ \tau + \mu_{(9,5,1)} := \dfrac{ \log{2}}{\log 10} -\dfrac{ \log{ ( 1+\dfrac{-4}{10})  (1/3) } } { {\log 10}} =1 \in \mathbb{Z},$$
and hence
\begin{eqnarray*}
\Gamma_{4} & = \left \lvert  d_1 \log 10 - (n-2) \log 2  -  \log (2/10) \right \rvert \\
& = \left \lvert  ( d_1 - 1) \log 10 - (n-3) \log 2  \right \rvert.
\end{eqnarray*}

If $ (a, b, {d_1-d_2 }   ) = (9,5,2) ,$ then
$$ 5\tau + \mu_{ (9,5,1) } := 5 \dfrac{ \log{2} }{ \log 10 } -\dfrac{ \log{ ( 1-\dfrac{4}{100} ) (1/3) } } {  \log 10 } =2 \in \mathbb{Z}.$$
In this case, we write
\begin{eqnarray*}
\Gamma_{4} & = \left \lvert  d_1 \log 10 - (n-2) \log 2  +  \log (32/100) \right \rvert \\
& = \left \lvert  ( d_1 - 2) \log 10 - (n-7) \log 2  \right \rvert,
\end{eqnarray*}

and hence
\begin{equation}
0<\left \lvert       \dfrac{\log 2 }{\log 10 }  - \dfrac{d_1 - 2}{n-7} \right \rvert  < \dfrac{4}{ 2^{k/2} {\log 10 } }.
\end{equation}

Similar to the first one,  we see that the bound  $ k<2000 $ is also valid in these two cases also.

We repeat the same reduction steps one more time but taking $ k<2000 $ and $ M:= 8.5 \cdot 10^{60} >n. $ When we work on \eqref{g3L}, this time, we take $ q_{129} $ instead of $ q_{588} $ and we find that
$$ \epsilon_a := ||\mu_a q_{129} || -M || \tau q_{129} || > 0.031955 ,$$ for each $ a \in \{1, 2, \ldots, 9 \} .$

By  Lemma \ref{reduction} we obtain  $ \lambda < 210 .$

Hence, $ \lambda=k/2 -5$ means that $ k <430.$ Assume that $ \lambda ={(d_1-d_2)} \dfrac{\log(10)}{\log(2)}-8 .$ Then
$$ {d_1-d_2} < 70.$$

Now, we pass to the $\Gamma_{4,}$ and we take $ q_{135} $ instead of $ q_{595} .$ Then we find that
$$\epsilon_{(a,b,d_1-d_2)}:= ||\mu_{(a,b,d_1-d_2)} q_{135} || -M || \tau q_{135} || \geq  \epsilon_{ (4,6,2) }  > 0.000065 $$
for each $ a \in \{1, 2, \ldots, 9 \} $, $ b \in \{1, 2, \ldots, 9 \} $, $ d_1-d_2 \in \{1, 2, \ldots, 70 \} $, except for the same three triples $(a, b, {d_1-d_2 })=(9,2,1), (9,5,1), (9,5,2) .$ Thus, we repeat the same calculations as we did before and we find that, even in the exceptional cases, $ k<460 $ which contradicts the fact $ k>470. $ So, we conclude that Equation \eqref{Lnab} has no solutions when $ k>470 $ and $ a \neq 0. $

\subsection{The Case $ a=0 $ and $k-$Lucas Numbers of the form $b10^{d_2}$}

Let $ a=0. $ Then \eqref{Lnab} turns into the equation
\begin{equation}\label{Lb10}
L_n^{(k)} =b10^{d_2}.
\end{equation}

Clearly, we take $ b > 0 .$ In fact, our previous work contains most of the material to solve this equation, with some small manipulation on the variables. So, in any applicable case, we follow the previous notation to prevent the recalculation. 

By \eqref{Lb10}, $ \Lambda_2 $ which was given in \eqref{lam2L} is  valid as
\begin{equation*}
\lvert \Lambda_2^{3} \rvert  :=  \lvert   \alpha^{-(n-1) } 10^{d_2} f_k(\alpha)^{-1} (2\alpha-1)^{-1} b -1     \rvert  \leq \dfrac{5}{2\alpha^{n-1}},
\end{equation*}
and $ \Lambda_2^{'} \neq 0. $ Let
$$ \eta_1:= \alpha, \eta_2:= 10, \eta_3:=f_k(\alpha)^{-1} (2\alpha-1)^{-1} b $$
with $ b_1:= -(n-1)$, $b_2:=d_2,$  $b_3:=1. $ So,
$$
h(\eta_3) \leq  h(f_k(\alpha)^{-1}) + h( (2\alpha-1)^{-1} ) +h(b)   \leq   3 \log{k} +\log 3+\log 9 < 8 \log k.
$$

From \eqref{Lb10} and \eqref{L1}, we may write $ 10^{d_2} \leq L_n^{(k)} \leq 2 \alpha^n < 2^n. $ Thus, it is enough to take $ B:= n-1. $ Note that, the inequalities $ 1+ \log k < 3 \log k$ and $ 1+ \log {n-1} < 2 \log {n-1}$ holds for all $ k \geq 2 $ and $ n \geq 4. $ We apply Theorem \ref{Matveev} by following the similar notation as we did before for $ \Lambda_2 ,$ we obtain that

$$ n-1 < 1.6 \times 10^{13} k^4 \log^2 k \log(n-1) .$$

We take $ T:=1.6 \times 10^{13} k^4 \log^2 k.$ Then $ \log T<60 \log k $ for all $ k \geq 2. $ Thus, from Lemma \ref{Guzman}, we find
\begin{equation}
\label{nfirst10}
n< 2.1 \times 10^{17} k^4 \log^4 k.
\end{equation}

Assume that $ k \leq 450 ,$ then $ n< 3 \times 10^{29}. $ By repeating the similar calculations, as we did before for \eqref{g2L} to the inequality,
\begin{equation*}
0 < \left \lvert (n-1) \dfrac{\log \alpha }{\log 10 }  -d_2 -  \dfrac{  \log( b f_k(\alpha)^{-1} (2\alpha-1)^{-1}   }{\log 10 }    \right \rvert < \dfrac{5}{\alpha^{n-1} \log 10  },
\end{equation*}
we see that the bounds found for $ a \neq 0 $ strictly hold for the case $ a=0 .$ Hence, by a computer search, we see that \eqref{Lb10} has no solution when  $ k \leq 470 .$  

Let $ k>450 .$ From \eqref{41L}, we write
$$ 0 \neq  \Lambda_4':= \left \lvert  2^{-(n-2)} 10^{d_2}b/3 -1 \right \rvert \leq  \dfrac{1}{2^{k/2}}  .$$

By taking
 
$$ (\eta_1, \lvert b_1 \rvert) :=(2, n-2),  (\eta_2, \lvert b_2 \rvert) :=(10, d_2) , (\eta_3, \lvert b_3 \rvert) :=(b/3, 1), $$
from Theorem~\ref{Matveev} together with \eqref{nfirst10}, we find $ k<4 \times 10^{14} $ and hence, from \eqref{nfirst10}, $ n<6.9 \times 10^{81}. $ To reduce these bounds, we write
$$ \Gamma_4' := \left \lvert (n-2) \log2 -d_2 \log 10 -\log (b/3) \right  \rvert , $$
so that, as we did before, we obtain
\begin{equation}
0< \left \lvert (n-2) \dfrac{ \log 2 }{\log 10}  - d_2 - \dfrac{ \log (b/3) }{ \log 10} \right \rvert < \dfrac{2}{ 2^{k/2} \log 10 }.
\end{equation}

Assume that $ b \not \in \{3,6 \}.$ 

Then, applying Lemma \ref{reduction} by choosing the parameters as $M:=6.9 \times 10^{81},$  $\mu_b:= - \log (b/3) /\log 10,$  
$\epsilon_{b}:= \lvert\lvert \mu_{b} q_{170} \rvert\rvert -M \lvert\lvert \tau q_{170} \rvert\rvert $ and the others as in the previous section, we find that $ k<564. $ If $ b $ is 3 or 6 then, from $\Gamma_4' $, we have that

$$ \left \lvert \dfrac{\log2}{\log10} - \dfrac{u}{v}  \right \rvert < \dfrac{2}{2^{k/2} v \log10}, $$
where $ \dfrac{u}{v} $ is $ \dfrac{d_2}{n-2} $ and $ \dfrac{d_2}{n-3} $, respectively. We use the theory of continued fractions as we did before for \eqref{gcf}, to obtain that $ k<572. $ Thus, from \eqref{nfirst10}, we obtain a reduced bound as $ n<4 \times 10^{31}. $ We repeat the same reduction algorithm with $ M:=4 \times 10^{31} $ and as a result we obtain that $ k<440,$ a contradiction. This completes the~proof.




{}

\end{document}